\documentclass[12pt]{article}
\usepackage[utf8]{inputenc}
\usepackage{amsmath}
\usepackage{amsfonts}
\usepackage{amssymb}
\usepackage{amsthm}
\usepackage{relsize}
\usepackage{array}
\usepackage{multirow}
\usepackage{tabu}
\usepackage{soul}
\usepackage{graphicx}
\usepackage{epigraph}
\usepackage{float}
\usepackage[english]{babel}
\usepackage[usenames, dvipsnames]{color}
\usepackage{fancyhdr}
\usepackage[Sonny]{fncychap}
\usepackage{tikz-cd}
\usepackage{tkz-euclide}
\usepackage{sseq}
\usepackage{newfloat}
\usepackage{amsmath,mathtools}
\usepackage{pgfplots}
\usepackage{geometry}
 \usepackage{sectsty}
\sectionfont{\fontsize{14}{15}\selectfont}

\usepackage[all]{xy}
\DeclareFloatingEnvironment[fileext=dia]{diagram}

\theoremstyle{definition}

\newtheorem{remark}{Remark}

\usepackage[most]{tcolorbox}

\makeatletter
\newtcolorbox{note}[1][]{%
	breakable,
	enhanced jigsaw, 
	borderline west={3pt}{0pt}{black!10!white}, 
	borderline south={1pt}{0pt}{black!10!white}, 
	borderline east={1pt}{0pt}{black!10!white},
	borderline north={1pt}{0pt}{black!10!white},
	sharp corners, 
	boxrule=0pt, 
	attach title to upper, 
	left=0pt,
	right=0pt,
	top=0pt,
	bottom=0pt,
	boxsep=5pt,
	colback=white,
	frame hidden,
	#1
}
\makeatother

\makeatletter
\newtcolorbox{note1}[1][]{%
	breakable,
	enhanced jigsaw, 
	sharp corners, 
	boxrule=0pt, 
	attach title to upper, 
	fontupper=\linespread{1.1}\fontfamily{qpl}\selectfont,
	fontlower=\linespread{1.1}\fontfamily{qpl}\selectfont, 
	left=0pt,
	right=0pt,
	top=0pt,
	bottom=0pt,
	boxsep=3pt,
	colback=green!3!white,
	frame hidden,
	before skip=10pt plus 2pt,after skip=10pt plus 2pt,
	#1
}
\makeatother

\makeatletter
\newcommand\tabfill[1]{%
	\dimen@\linewidth
	\advance\dimen@\@totalleftmargin
	\advance\dimen@-\dimen\@curtab
	\parbox[t]\dimen@{#1\ifhmode\strut\fi}%
}
\makeatother

\usepackage{tabularx}
\usepackage{imakeidx}
\usepackage{hyperref}
\hypersetup{colorlinks={true},linkcolor={blue},citecolor={blue},urlcolor={blue}, linktoc=all}  
 
 \usepackage[capitalize, nameinlink]{cleveref}
 \crefdefaultlabelformat{#2\textbf{#1}#3}
 \crefname{figure}{Figure}{Figures} 
 \Crefname{figure}{Figure}{Figures}
 \crefname{table}{Table}{Tables}
 \Crefname{table}{Table}{Tables}
 \crefname{section}{\S\hspace{-1mm}}{\S\hspace{-1mm}}
 \Crefname{section}{\S\hspace{-1mm}}{\S\hspace{-1mm}}
 \crefname{equation}{}{}
 \Crefname{equation}{}{}
 \crefname{example}{Geometric Pattern}{Geometric Patterns} 
 \Crefname{example}{Geometric Pattern}{Geometric Patterns}

 \usepackage{graphicx,tipa}

 \usepackage{footnote}

\begin{document}

\title{\textbf{Pythagorean Theorem in Elamite Mathematics}}

\author{Nasser Heydari\footnote{Email: nasser.heydari@mun.ca}~ and  Kazuo Muroi\footnote{Email: edubakazuo@ac.auone-net.jp}}

\maketitle

\begin{abstract}
This article studies the application of the   Pythagorean theorem in the Susa Mathematical Texts (\textbf{SMT}) and we discuss those texts whose problems and related calculations demonstrate its use. Among these texts,   \textbf{SMT No.\,1} might be the most important as it contains a geometric application of the  Pythagorean theorem.  
\end{abstract}

\section{Introduction}
The Pythagorean theorem is used throughout the \textbf{SMT}--both explicitly and implicitly. Here, we  consider  only those applications of the theorem found in  \textbf{SMT No.\,1}, \textbf{SMT No.\,3}, \textbf{SMT No.\,15} and \textbf{SMT No.\,19}. These texts were inscribed by Elamite scribes between 1894--1595 BC on   26 clay tablets excavated from Susa in  southwest Iran by French archaeologists in 1933. The texts of all the tablets,  along with their interpretations, were first published in 1961 (see \cite{BR61}).

On the obverse of \textbf{SMT No.\,1}\footnote{The reader can see   this tablet on the website of the Louvre's collection. Please see \url{https://collections.louvre.fr/en/ark:/53355/cl010185651} for obverse  and   reverse.}, there is an isosceles triangle  inscribed in a circle  along with  some numerical data. The reverse has several fragmentary numbers whose relations to the figure on the obverse unfortunately are not clear.

The text of  \textbf{SMT No.\,3}\footnote{The reader can see   this tablet on the website of the Louvre's collection. Please see \url{https://collections.louvre.fr/en/ark:/53355/cl010185653} for obverse  and   reverse.}    comprises a list of constants with 68 entries, almost all of which are preserved in good condition. The content of this tablet is as follows: line 1 is the headline, lines 2-32 contain the mathematical constants regarding areas and dimensions of geometrical figures  and  lines 33-69 contain the non-mathematical constants  concerning work quotas, metals, and so on. The last two lines (lines 70-71)  are an example of how to use the constants  of work quotas.

\textbf{SMT No.\,15}\footnote{The reader can see   this tablet on the website of the Louvre's collection. Please see \url{https://collections.louvre.fr/en/ark:/53355/cl010186541} for obverse  and   reverse.} contains three problems, two similar ones on the obverse and a badly damaged one on the reverse. Although the first and second problems are comparatively well preserved, they remain  unintelligible to us, as do their solutions to date. It is clear, however, that they are applied problems of the Pythagorean Theorem concerned with the enlargement of a gate. 
 
The text  of \textbf{SMT No.\,19}\footnote{The reader can see   this tablet on the website of the Louvre's collection. Please see \url{https://collections.louvre.fr/en/ark:/53355/cl010186429} for obverse  and   reverse.} contains two problems, one on the obverse and the other on the reverse of the tablet, both of which deal with simultaneous equations concerning Pythagorean triples. In the first problem the diagonal (of a rectangle) or the hypotenuse (of a right triangle) is called {\fontfamily{qpl}\selectfont tab}  ``friend, partner'' which reminds us of the fact that a Pythagorean triple was called {\fontfamily{qpl}\selectfont illat} ``group, clan'' in Babylonian mathematics. Moreover,   the scribe of this tablet handles these equations skillfully, especially in the second problem which might  be one of the most complicated systems of simultaneous equations in  Babylonian mathematics.

\section{Pythagorean Theorem}

\subsection{Statement}
\begin{figure}[H]
	\centering
	\includegraphics[scale=1]{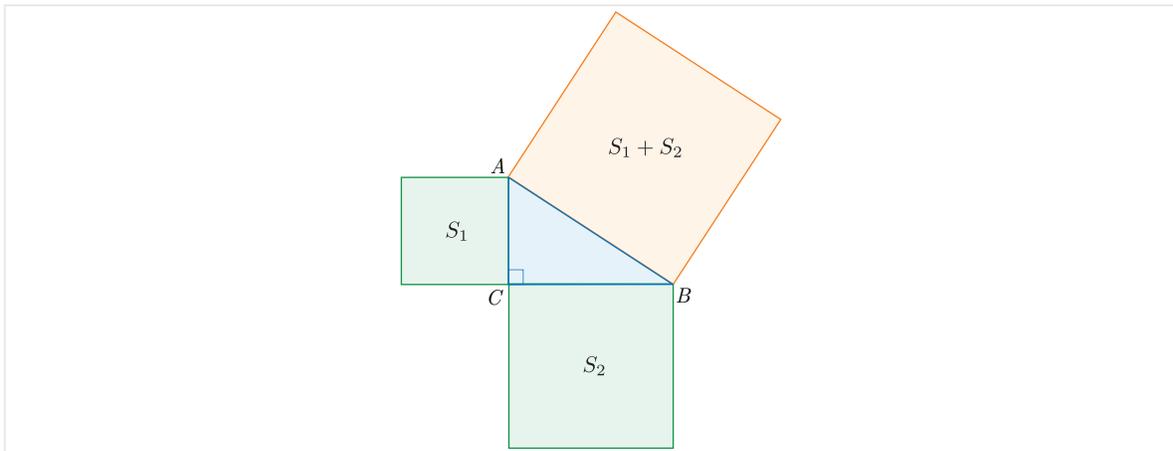}
	\caption{Pythagorean theorem}
	\label{Figure1}
\end{figure}

One of the oldest elementary mathematical theorems taught to students in middle or high school is the \textit{Pythagorean theorem}. Geometrically, this theorem states that for any given right triangle, the area of the square whose side is the hypotenuse (that opposite the right angle) is equal to the sum of the areas of the two squares whose sides form the other two legs of the right triangle.  A pictorial representation of this theorem is given in \cref{Figure1} in which the total   area  of the two smaller green squares  is equal to  the area of the  bigger orange one.

The best-known  statement of this theorem is usually given as an algebraic equation with respect to the lengths of sides of a right triangle.  Consider a right triangle $ \triangle ABC $ such that its angle $ \angle ACB=90^{\circ}$ and  the lengths of its two legs are given by $\overline{BC}=a $, $\overline{AC}=b $ and that of its hypotenuse   by $\overline{AB}=c $.   Then the Pythagorean theorem simply states  that 
\begin{equation}\label{equ-SMT1-a-a}
	c^2=a^2+b^2.	
\end{equation}
This is usually called the \textit{Pythagorean rule} or \textit{Pythagorean formula}. It should be noted that the converse of the Pythagorean theorem is also true in that if the lengths of three sides of a triangle satisfy \cref{equ-SMT1-a-a}, then the triangle has a right angle.

\subsection{History}
\begin{figure}[H]
 	\centering
 	\includegraphics[scale=1]{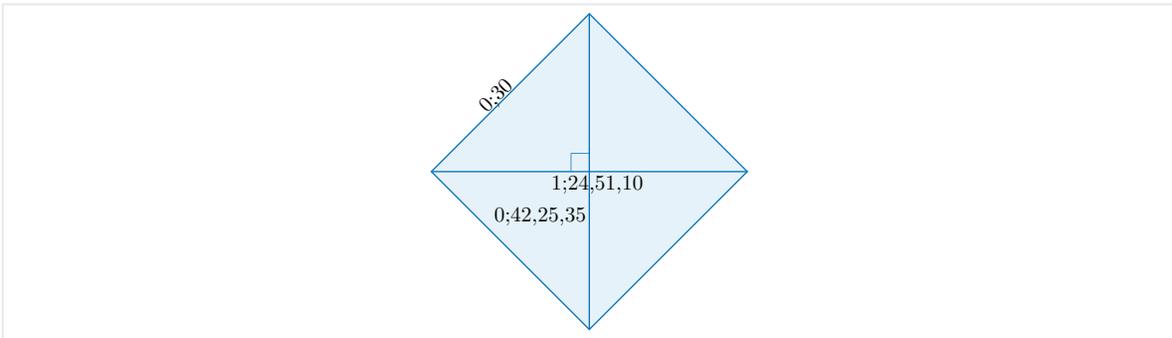}
 	\caption{Reconstruction of \textbf{YBC 7289}}
 	\label{Figure3}
 \end{figure}

 Although this theorem is named after the  famous Greek philosopher and mathematician Pythagoras  of Samos (circa 570--495 BC), its origin  dates   to millennia before him\footnote{For a discussion on the origin of the Pythagorean theorem in the Babylonian mathematics, see \cite{Hyp99}.}. It is  well-known   that the Babylonian and Elamite scribes   were familiar with this theorem long before the Greeks, and  there are  a number of   their clay tablets containing applications of this theorem\footnote{For a list of   known   Babylonian applications of the Pythagorean theorem, see \cite{Fri07-1}.}. For example, in the   tablet  \textbf{YBC 7289}\footnote{The mathematical tablet \textbf{YBC 7289}  belongs to the Yale Babylonian Collection. Its text and interpretation were published by Neugebauer in \cite{NS45}. For photos of its obverse and reverse, see \url{https://commons.wikimedia.org/wiki/File:YBC-7289-OBV-REV.jpg}.}, there is a square with the numbers written on its sides and diagonals  shown in \cref{Figure3}.

It has been suggested that the scribe of this tablet has used the relation $b=a\sqrt{2}$ and an approximation to  $ \sqrt{2}$ for finding the diagonal $b$ of the square with side $a$ (see \cite{FR98,NS45,Neu69}  for more details). It is clear that the relation $b=a\sqrt{2}$ follows from the Pythagorean theorem in an isosceles right triangle with sides $a,a,b $. In fact, by using the Pythagorean formula for the sides of the upper right triangle in \cref{Figure3}, we have  $b^2=a^2+a^2=2a^2$  implying that   $b=a\sqrt{2}$

\begin{figure}[H]
	\centering
	\includegraphics[scale=1]{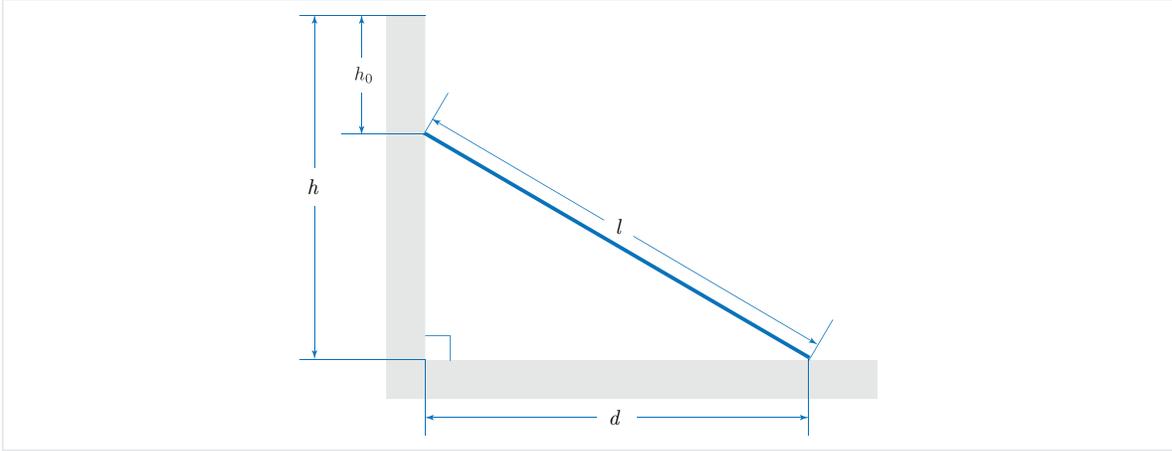}
	\caption{A timber against a wall}
	\label{Figure4}
\end{figure}

 Another example is given  in a problem from the mathematical  tablet \textbf{BM 85196}\footnote{This mathematical tablet, which is held in the British Museum contains 18 problems on a variety of subjects. The text of this tablet was originally published by Thureau-Dangin. For more information about its text, see \cite{Hyp02, NS45,Thu35}.} in which a timber of length 30 stands against a wall of height 30 such that the upper end has slipped down by 6. In this text, the scribe   computes  the distance  the lower end moves using the Pythagorean formula\footnote{For a more detailed discussion on this problem, see \cite{Mur91-2} or \cref{appendix}  of this article.}.  \cref{Figure4} depicts such a situation and by using the Pythagorean theorem for the right triangle formed by the timber, the wall and the ground, one  gets
\begin{equation}\label{equ-SMT1-a-b}
	l^2=d^2+(h-h_0)^2.
\end{equation}

  In most cases of such problems, the values of $l$, $h$ and $h_0$ are known and the scribe uses   equation \cref{equ-SMT1-a-b} to find  the value of $d$:
\begin{equation*}
	d=\sqrt{l^2-(h-h_0)^2}.
\end{equation*}

\subsection{Applications}
The Pythagorean theorem has different applications in mathematics. Here, we only list some of the most important applications.

\subsubsection*{Construction of Incommensurable Lengths}
One of the earliest applications of this theorem is to construct  lengths like $\sqrt{2}, \sqrt{3}, \sqrt{5}$ and so on.  Two non-zero real numbers $a$ and $b$ are called \textit{incommensurable} if their ratio is not a rational number. For example, any pair $(\sqrt{n},1)$, where $n$ is not a perfect square, is incommensurable. 

\begin{figure}[H]
	\centering
	\includegraphics[scale=1]{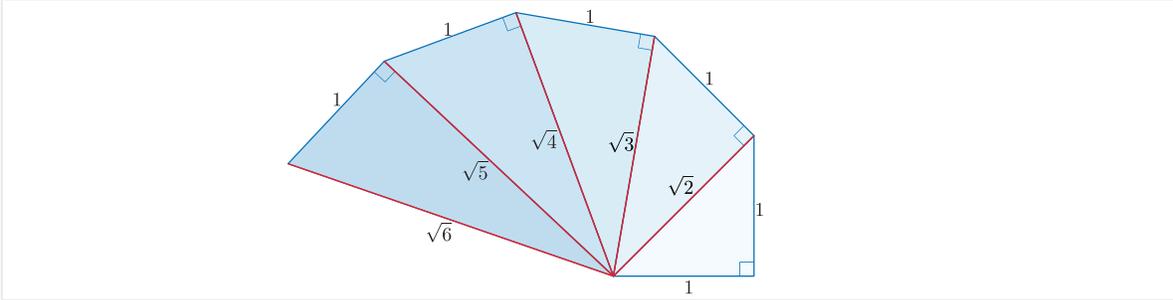}
	\caption{Incommensurable lengths}
	\label{Figure8}
\end{figure}

Note that for any natural number $n$, we always can write  
\[ \sqrt{n}=\sqrt{n-1+1}=\sqrt{(\sqrt{n-1})^2+1^2}. \]
This says that the length $\sqrt{n}$ is constructible by straightedge and compass, if we  can construct $\sqrt{n-1}$. In that case, all we need to do is to construct a right triangle with legs $1$ and $\sqrt{n-1}$. The hypotenuse is then the square root $\sqrt{n}$ (see \cref{Figure8}).

\subsubsection*{Euclidean Distance}
Perhaps one of  the main applications of this  theorem appears in   Euclidean geometry where it serves as a basis for the definition of the distance between two points in the plane. The distance $ d(A,B)$ between two points $A=(x_1,y_1)$ and $B=(x_2,y_2)$ in the $xy$-coordinate system is defined as
\begin{equation}\label{equ-SMT1-a-aa}
	d(A,B)=\sqrt{(x_2-x_1)^2+(y_2-y_1)^2}.	
\end{equation}
Note that a similar formula is also used for Euclidean spaces of higher dimensions. 

\begin{figure}[H]
	\centering
	\includegraphics[scale=1]{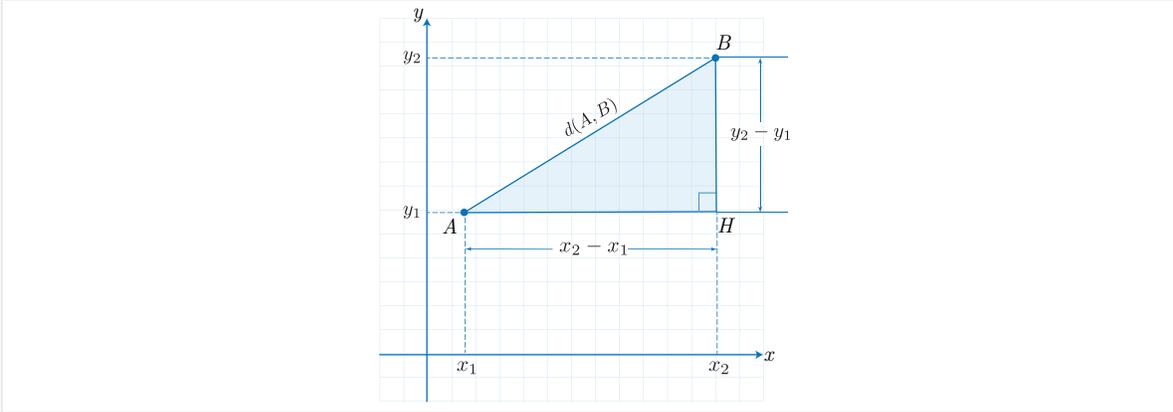}
	\caption{Euclidean distance}
	\label{Figure2}
\end{figure}

It is clear from \cref{Figure2} that the distance $d$ is the length of the hypotenuse $ AB$ in the blue right triangle $\triangle AHB $ whose legs are of lengths $\overline{AH}=x_2-x_1$ and $\overline{BH}=y_2-y_1$. Formula \cref{equ-SMT1-a-aa} is a direct consequence of formula \cref{equ-SMT1-a-a}.

\subsubsection*{Trigonometry}
\begin{figure}[H]
	\centering
	\includegraphics[scale=1]{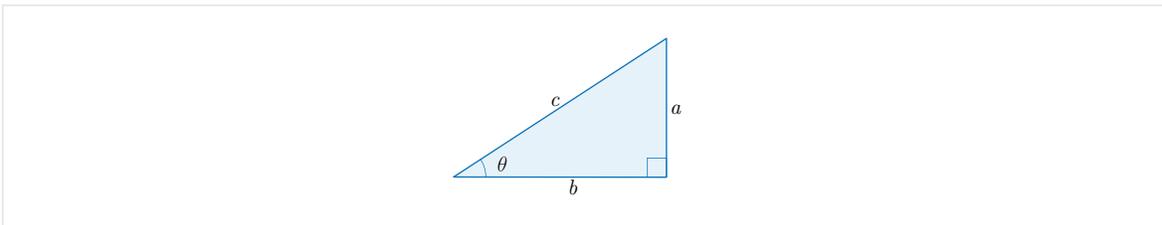}
	\caption{Trigonometric relations}
	\label{Figure6}
\end{figure}

Trigonometry is based on right  triangles and the relation between their sides and angles. The two main trigonometric relations, i.e., the  \textit{sine} and the \textit{cosine} of an acute angle $\theta$ in a right triangle with legs $a,b$ and hypotenuse $c$ are defined as $\sin(\theta)=\frac{a}{c}$ and $\cos(\theta)=\frac{b}{c}$ (see \cref{Figure6}). The Pythagorean theorem simply implies the most famous trigonometric identity:
\[ \sin^2(\theta)+\cos^2(\theta)=1. \]
In a sense, this identity and the Pythagorean rule are equivalent.

\subsubsection*{Complex Numbers}
\begin{figure}[H]
	\centering
	\includegraphics[scale=1]{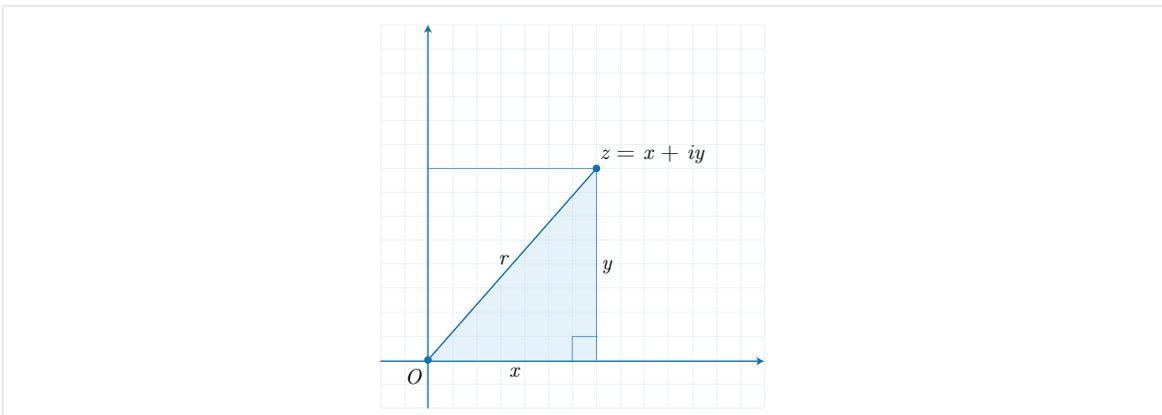}
	\caption{Complex numbers}
	\label{Figure7}
\end{figure}

Complex numbers are defined as $z=x+iy$ where $x,y$ are real numbers and $i=\sqrt{-1}$ is the \textit{imaginary unit}. Any complex number   $z=x+iy$ can be identified with the point $(x,y)$ in the two-dimensional coordinate system (see \cref{Figure7}). In that case, the distance between the point and the origin is called the \textit{absolute value} of $z$ and denoted by $r$. By the Pythagorean theorem, we always have $r^2=x^2+y^2$.

\subsection{Pythagorean Triples}
Besides the geometric aspects of the Pythagorean formula which concern   right triangles, we can also consider it from a purely algebraic point of view. In fact, any triple $(a,b,c)$ of  positive integers satisfying formula \cref{equ-SMT1-a-a} is called a \textit{Pythagorean triple}. Note that if $(a,b,c)$ is a Pythagorean triple, for any natural number $n>1$, the new triple $(na,nb,nc) $ is also Pythagorean, because it clearly satisfies formula \cref{equ-SMT1-a-a}.   If three integers in a Pythagorean triple $(a,b,c)$ do not have a common factor, it is called a \textit{primitive} Pythagorean triple. For example,    triples such as $(3,4,5)$ and $(5, 12, 13) $ are primitive Pythagorean triples while $(6,8,10)=(2\times 3, 2\times 4, 2\times 5)$ is not primitive. The following is a list of some primitive Pythagorean triples:

\begin{center}
	\begin{tabular}{lllll}
		(3,4,5)&	(5,12,13)&	(7,24,25)&	(8,15,17)&	(9,40,41)\\
		(11,60,61)&	(12,35,37)&	(13,84,85)&	(15,112,113)&	(16,63,65)\\
		(17,144,145)&	(19,180,181)&	(20,21,29)&	(20,99,101)&	(21,220,221)\\
		(23,264,265)&	(24,143,145)&	(25,312,313)&	(27,364,365)&	(28,45,53)\\
		(28,195,197)&	(29,420,421)&	(31,480,481)&	(32,255,257)&	(33,56,65)\\
		(33,544,545)&	(35,612,613)&	(36,77,85)&	(36,323,325)&	(37,684,685)\\
		(39,80,89)&	(39,760,761)&	(40,399,401)&	(41,840,841)&	(43,924,925)
	\end{tabular}
\end{center}

One of the classical problems in number theory was to determine  parametric   formulas generating primitive Pythagorean triples for different values of   the parameter. A fundamental formula was provided by    Euclid of Alexandria  (circa 325--265 BC)  saying that if $m>n>0$ are two coprime\footnote{Two natural numbers are called coprime if their only common factor is 1.} integers both of which are not odd, then the triple $$(m^2-n^2, 2mn,  m^2+n^2)$$
provides us with a primitive Pythagorean triple. Note that
\[ (m^2+n^2)^2=(m^2-n^2)^2+(2mn)^2. \]

The inverse of this statement is also true: for any primitive Pythagorean triple $(a,b,c)$, there exist positive integers $m,n$ satisfying the above-mentioned conditions such that $a=m^2-n^2$, $b=2mn$, and $c=m^2+n^2$.  Although this formula only gives  primitive Pythagorean triples,  the modified formula $(km^2-kn^2, 2kmn,  km^2+kn^2)$, where $k$ is a natural number, provides all  the Pythagorean triples uniquely. One can also consider the Pythagorean triples in a purely geometric point of view. In fact, any such triple $(a,b,c)$ can be identified with the integer point $(a,b)$ in the coordinate plane such that its distance to the origin is a positive integer, because $c=\sqrt{a^{2}+b^2}$.

Similar to the Pythagorean theorem, the Pythagorean triples were known to the ancient mathematicians too. One of the most famous Babylonian mathematical clay tablets, \textbf{Plimpton 322}\footnote{This mathematical tablet  is one of the most famous Babylonian clay tablets whose cuneiform text was first published by Neugebauer  in \cite{NS45}. It has a table of four columns and 15 rows of numbers which most scholars believe to be  Pythagorean triples. For a detailed discussion about the text of this tablet, see \textit{Babylonian Number Theory and Trigonometric Functions: Trigonometric Table and Pythagorean Triples in the Mathematical Tablet Plimpton 322} by K. Muroi published in \cite{KKL13}, pages 31-47.}, which has been under discussion for many years, is widely believed to contain  a list of  fifteen Pythagorean triples.  Besides this text,   in the text of \textbf{SMT No.\,19} the  Susa scribes deal  with the  two Pythagorean triples $(24,32,40)$ and $(30,40,50)$ which are obtained   from the primitive Pythagorean triple $(3,4,5) $.   The primitive Pythagorean triple $(7,24,25)$   also occurs   in the mathematical calculations regarding  the text of \textbf{SMT No.\,1} and  \textbf{SMT No.\,3}.

\subsection{Proof}
The proof of the Pythagorean formula \cref{equ-SMT1-a-a} has been of great interest to  many mathematicians through the ages. It is believed that the first proof was   provided by   Euclid of Alexandria and from that time forth     many    people (mathematicians and non-mathematicians such as the famous Leonard Da Vinci  and the 12-years old Einstein)   have sought  to provide new proofs for this theorem. As is stated in  \cite{Mao07}, there are at least four hundred proofs for this   formula using different  approaches some of which are set out in that book.\footnote{Besides  a few proofs   given in \cite{Mao07},  on the website \url{https://www.cut-the-knot.org/pythagoras/index.shtml}, the reader can find 121 proofs for the Pythagorean formula. For a full list of 371 proofs, see \cite{Loo72}.}

\begin{figure}[H]
	\centering
	\includegraphics[scale=1]{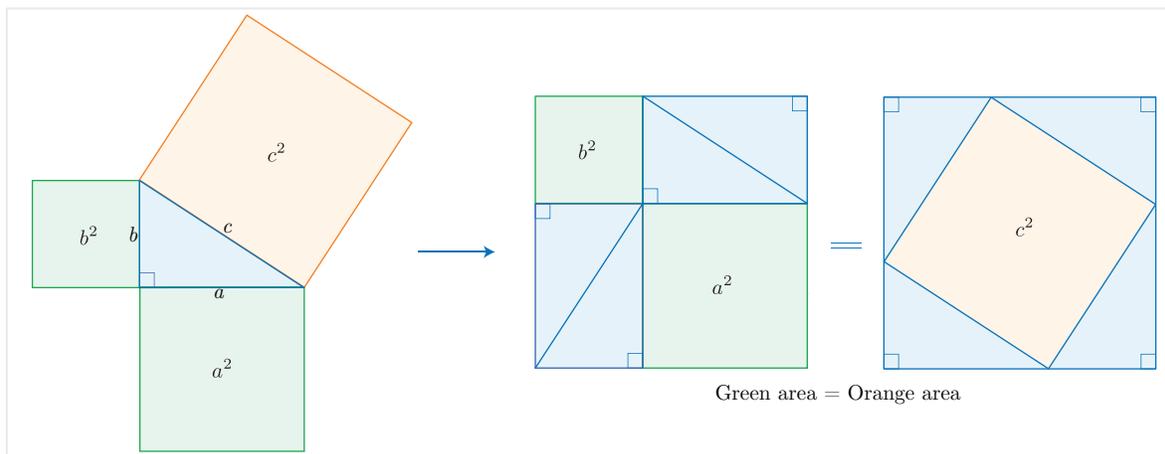}
	\caption{A geometric proof for the Pythagorean theorem}
	\label{Figure5}
\end{figure}

There are many  proofs  using elegant geometric  reasoning. One of these   is shown in \cref{Figure5}. As is seen from the figure, we can form a square of side $a+b$ with four copies of the right triangle $a,b,c$ and two squares of sides $a$ and $b$. At the same time, one can also form the same square of side $a+b$  with the four copies of the right triangle $a,b,c$ and a square of side $c$. If we remove the four right triangles in each layout, the remaining parts are equal, meaning that $c^2=a^2+b^2$.

\subsection{Generalization}
\begin{figure}[H]
	\centering
	\includegraphics[scale=1]{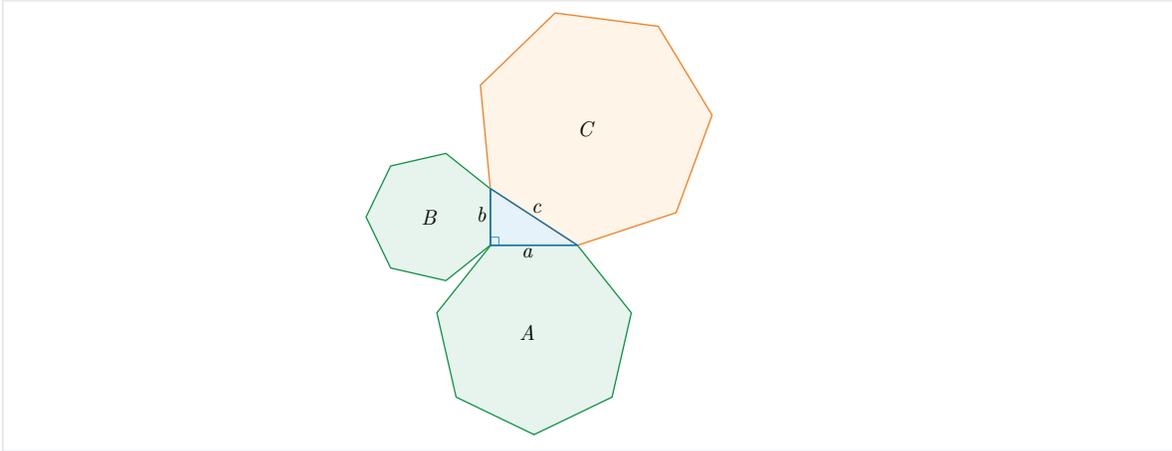}
	\caption{A generalization of the Pythagorean theorem: $A+B=C$}
	\label{Figure5-1}
\end{figure}

It is interesting that if we replace the squares  in the Pythagorean theorem with similar shapes, the relation between their areas still holds. In other words, if we  build similar figures   on the sides of a right triangle, then the area of the figure on the hypotenuse is equal to the sum of those of the two other figures on the legs. \cref{Figure5-1} shows this situation for similar regular heptagons built on the sides of a right triangle.

\section{Pythagorean Theorem in the SMT}
In this section, we discuss the applications  of the Pythagorean theorem in the Susa Mathematical Texts (\textbf{SMT}). In some cases, the scribe has implicitly used this theorem, but there are explicit expressions of the Pythagorean rule in some texts. We have only considered the most important examples here, although the traces of this basic theorem can be seen in many texts of the \textbf{SMT}  (see \cite{HM22-1,HM22-2,HM22-3,HM23-1,HM23-2}).

\subsection{\textbf{SMT No.\,1}}\label{SMT1}

\subsubsection*{Transliteration}

\begin{note1} 
	\underline{Obverse: Lines 1-6}
	\vspace{-3mm}
	\begin{tabbing}
		\hspace{10mm} \= \kill   
		(L1)\> \tabfill{50 u\v{s}}\\
		(L2)\> \tabfill{bar-d\'{a}}\\
		(L3)\> \tabfill{30}\\
		(L4)\> \tabfill{31;15  u\v{s}}\\
		(L5)\> \tabfill{8;45}\\
		(L6)\> \tabfill{[40]  u\v{s} sag-bi \textit{ga-am-ru}}
	\end{tabbing}
\end{note1}

\noindent
According to the reverse, the numbers we can recognize are as follows:

\begin{note1} 
	\underline{Reverse}  
	$$ \text{1,33,27,30\hspace{7mm}  [$\cdots$]40(?),8,45\hspace{7mm}   [5]2,30\hspace{7mm}  8,14}$$ 
\end{note1}

\subsubsection*{Translation}

\underline{Obverse: Lines 1-6}
\vspace{-3mm}
\begin{tabbing}
	\hspace{10mm} \= \kill   
	(L1)\> \tabfill{``$50$ is the length'': written over $AC$ and under $BC$.}\\
	(L2)\> \tabfill{``the transversal ($AB$)'': written over $BM$. The word {\fontfamily{qpl}\selectfont bar-d\'{a}} is a variant of the  Sumerian word {\fontfamily{qpl}\selectfont bar-da} ``crosspiece, crossbar''. In the \textbf{SMT}, it is used in the  sense of ``transversal,  diagonal''.}\\
	(L3)\> \tabfill{$30$: written (as the length of $BM$) under $BM$.}\\
	(L4)\> \tabfill{``$31;15$ is the length'': written under $BO$.}\\
	(L5)\> \tabfill{$8;45$ ($=40-31;15$): written under $MO$.}\\
	(L6)\> \tabfill{``$40$ is the complete length of the vertex (that is $MC$)'': written  under $OC$.}
\end{tabbing}

\noindent 
\underline{Reverse}  
\begin{align*}
	1,33,27,30\\
	[\cdots]40(?),8,45\\
	[5]2,30\\
	8,14	
\end{align*}

\noindent
Regrettably, we do not as yet understand how these numbers relate to the figure on the obverse.

\subsubsection*{Mathematical Interpretation}
\begin{figure}[H]
	\centering
	\includegraphics[scale=1]{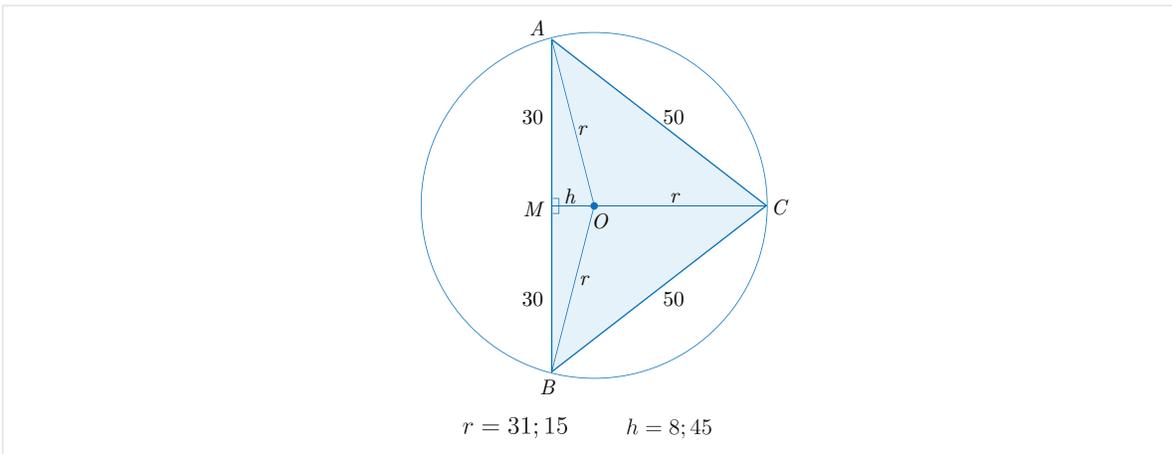}
	\caption{An isosceles triangle inscribed in a circle}
	\label{Figure11}
\end{figure}

We have reconstructed the drawing on the obverse of \textbf{SMT No.\,1}  in \cref{Figure11}. Here, $\overline{CA}=\overline{CB}=50$ and $\overline{AM}=\overline{BM}=30$. It seems that the scribe has intended to find the values of the radius $r$ and height $h$. There are different ways to solve this problem and find the values of $r$ and $h$.

 First, to compute the   height  $ \overline{CM} $ of the isosceles triangle  $\triangle ABC $, one can   use  the Pythagorean theorem   in the right triangle  $\triangle AMC$:
\begin{align*}
	\overline{CM}&=\sqrt{\overline{CA}^2-\overline{AM}^2}\\
	&=\sqrt{50^2-30^2}\\
	&=\sqrt{41,40-15,0}\\
	&=\sqrt{26,40}\\
	&=\sqrt{40^2}\\
	&=40.
\end{align*}
So,
\begin{equation}\label{equ-SMT1-a}
	\overline{CM}=40,
\end{equation}
and since $ h=\overline{OM}=\overline{CM}-\overline{OC}$, we also get
\begin{equation}\label{equ-SMT1-aa}
	h= 40-r.
\end{equation}
Next, to find the value of radius $r$,   we can   use the Pythagorean theorem  in the right triangle  $\triangle BMO$ to solve an equation with respect to $r$ as follows:
\begin{align*}
	~~&~~\overline{OB}^2= \overline{OM}^2+\overline{MB}^2\\
	\Longrightarrow~~&~~r^2 =(40-r)^2+30^2    \\
	\Longrightarrow~~&~~r^2 =40^2-(2\times 40) r+r^2+30^2  \\
	\Longrightarrow~~&~~r^2 =41,40-(1,20) r+r^2  \\
	\Longrightarrow~~&~~(1,20) r=41,40  \\
	\Longrightarrow~~&~~r=\frac{1}{(1,20)} \times (41,40)  \\
	\Longrightarrow~~&~~r=\left(0;0,45\right) \times (41,40)  \\
	\Longrightarrow~~&~~r=31;15. 
\end{align*}
Therefore, we get the value of the radius of the circumscribed circle 
\begin{equation}\label{equ-SMT1-b}
	r=31;15.
\end{equation}
Finally, it  immediately follows  from \cref{equ-SMT1-aa} and \cref{equ-SMT1-b} that 
\begin{align*}
	h&=40-r\\
	&=40-31;15\\
	&=8;45.
\end{align*}

\begin{remark}\label{rem-SMT1-a}
	Note that there are two primitive Pythagorean triples hidden in   this text: $ (3,4,5) $ and  $ (7,24,25) $. In fact,  in the right triangle $ \bigtriangleup AMC $,  the triple is 
	\[ (30,40,50)=(10\times 3,10\times 4,10\times 5) \]
	while in the right triangle $ \bigtriangleup AMO $,  the triple is 
	\[ \left(8\frac{3}{4},30,31\frac{1}{4}\right)=\left(\frac{5}{4}\times 7,\frac{5}{4}\times 24,\frac{5}{4}\times 25\right). \]
\end{remark}

\begin{remark}\label{rem-SMT1-b}
	Similar interpretations have been given by H\o{}yrup and Friberg in \cite{Hyp02,Fri07-1,Fri07-2}.
\end{remark}

\subsection{\textbf{SMT No.\,3}}
Although   the Pythagorean theorem is necessary in the calculation of many of geometric constants in \textbf{SMT No.\,3}  (see \cite{HM22-1,HM22-3}), we have only considered the main ones here.

\subsubsection*{Line 29: Height of an Equilateral Triangle}
In this line, we read {\fontfamily{qpl}\selectfont 52,30 igi-gub \textit{\v{s}\`{a}} sag-d\`{u}} ``0;52,30 is the constant of a equilateral triangle''. This number is the height of an equilateral triangle of side 1. 

\begin{figure}[H]
	\centering
	\includegraphics[scale=1]{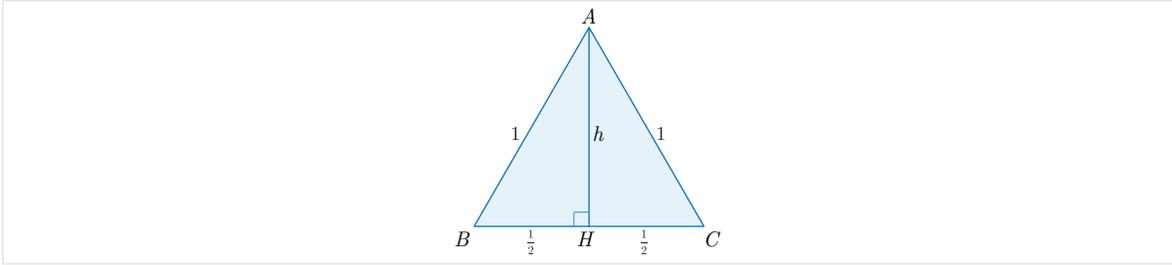}
	\caption{Height of an equilateral triangle}
	\label{Figure9}
\end{figure}

In fact, since the height $h$ bisects the base, in the equilateral triangle with side 1, we can use the Pythagorean theorem for a right triangle whose hypotenuse is 1 and whose base is $\frac{1}{2}$:
\[ h^2+\left(\frac{1}{2}\right)^2=1^2 \Longrightarrow h=\sqrt{1-\frac{1}{4}}=\frac{\sqrt{3}}{2}. \]
If we use the common Babylonian approximation $\sqrt{3}\approx \frac{7}{4}$, we get
\[ h\approx \frac{\frac{7}{4}}{2} = \frac{7}{8}=0;52,30\]
confirming the scribe's value.

\subsubsection*{Line 30: A right Triangle}
In this line, we read {\fontfamily{qpl}\selectfont 57,36 igi-gub \textit{\v{s}\`{a}} ub} ``0;57,36 is the constant of a right triangle''.  We have interpreted the number  $0;57,36 = \frac{7}{25}$  to be  one  side of a right triangle  whose  hypotenuse is     1 and whose other side is  $\frac{24}{25}$, because the Pythagorean rule holds in this triangle:
\begin{align*}
	\left(\frac{7}{25}\right)^2+\left(\frac{24}{25}\right)^2   =\frac{7^2+24^2}{25^2}  =\frac{25^2}{25^2}=1. 
\end{align*}

Note that these three numbers are multiplications of the Pythagorean triple $(7,24,25)$. A right triangle  similar to this one, whose sides have lengths  $a=31;15=\frac{125}{4}$,   $b=8;45=\frac{35}{4}$ and  $c=30 $ also occurs in \textbf{SMT No.\,1} (see \cref{SMT1}). Note that these numbers are multiplications of the Pythagorean triple  $(7,24,25)$:
\[ \frac{35}{4}=\frac{5}{4}\times 7,~~~30= \frac{5}{4}\times 24,~~~\text{and}~~~\frac{125}{4} =\frac{5}{4}\times 25. \]

\begin{figure}[H]
	\centering
	\includegraphics[scale=1]{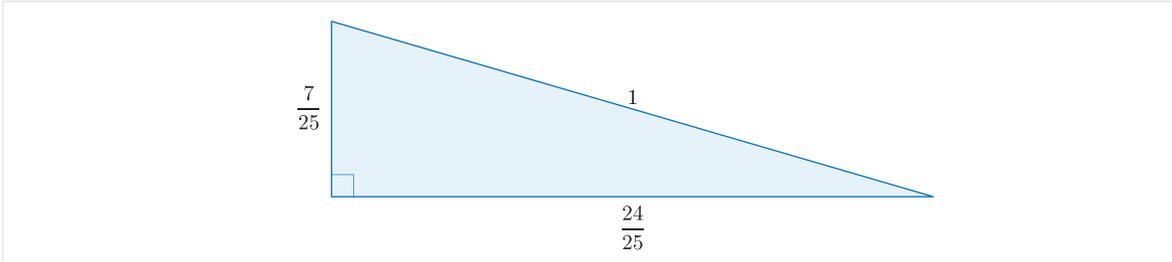}
	\caption{A right triangle}
	\label{Figure10}
\end{figure}   

\noindent
\textbf{Controversial number $0;57,36 $ and  approximation   $\pi\approx 3\frac{1}{8}$}\\
The number $\frac{24}{25} =0;57,36 $ in  line 30 is probably one of the most disputed constants in Babylonian mathematics  because    many interpretations  have been provided by  different scholars (see \cite{Bru50,BR61,Fri07-1,Mur92-1,Neu69}, for example). The controversy was primarily sparked off by Bruins in \cite{Bru50,BR61}, where he interpreted this number as a constant or the area of a circle. This resulted in the approximate value $\pi\approx \frac{25}{8}=3.125 $  for $\pi$ which is   more accurate than   the common Babylonian approximate value $\pi\approx 3$  and  even  the famous Egyptian approximate vale   $\pi\approx \frac{256}{81}$.

Although     disputes arose between scholars over this interpretation of Bruins,  it is interesting  that in  many books on the history of $\pi$  (see \cite{AH01,Bec71,PL04}, for example), this number has been referred to as one of the first approximations of $\pi$  in   ancient times (for a discussion about this alleged value for $\pi$, please see \cite{Mur92-1}). Whereas there are doubts about attributing this   approximate  value $\pi=3.125=3;7,30$ in \textbf{SMT No.\,3}  to Susa scribes,   it seems that Sumerian scribes  knew this approximate value. In fact,  there is another fragmentary tablet published   by Thureau-Dangin  \cite{Thu03}   on which a circle  has been drawn and  the  number 7 {\fontfamily{qpl}\selectfont sar} 10 {\fontfamily{qpl}\selectfont gín}    has been written  as the area  of  a circular plot.  

\begin{figure}[H]
	\centering
	\includegraphics[scale=1]{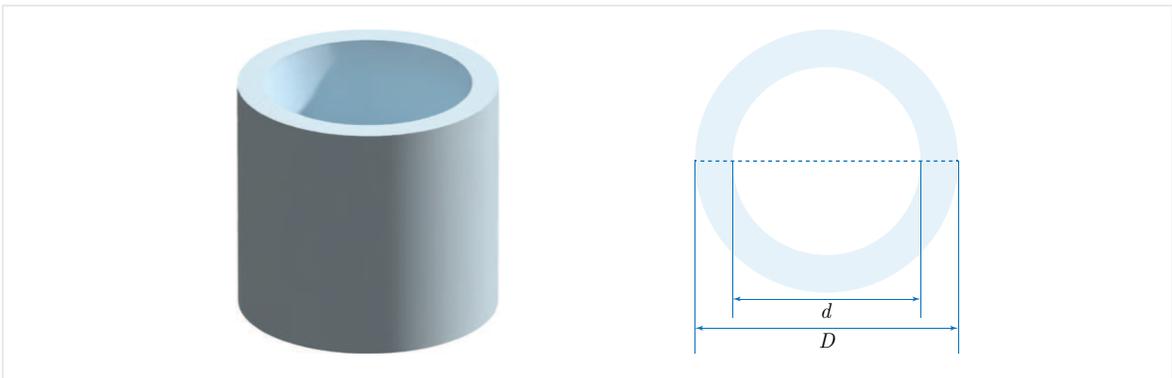}
	\caption{Top view of a cylindrical storehouse}
	\label{Figure10a}
\end{figure}

   An interpretation in \cite{Mur16} offers that the scribe of this tablet has used the approximation $ \pi=3;7,30$. In fact, if   the inner diameter of the storehouse  is $d=3$ {\fontfamily{qpl}\selectfont nindan} and the thickness of its wall is $\frac{1}{2}$ {\fontfamily{qpl}\selectfont šudùa}\footnote{Note that 36 {\fontfamily{qpl}\selectfont šudùa} is equal to 1 {\fontfamily{qpl}\selectfont nindan}}, then its external diameter is $D=3+2\times (0;0,50)=3;1,40$ {\fontfamily{qpl}\selectfont nindan}. It follows from $S=\frac{c^2}{4\pi}$ and $c=\pi D$  that $S=\frac{\pi D^2}{4} $, so 
\begin{align*}
	\pi&=\frac{4S}{D^2}\\
	&\approx\frac{4\times (7;10)}{\left(3;1,40\right)^2}\\ 
	&\approx (28;40)  \times \frac{1}{(9;10,2,46,40)}\\ 
	&\approx (28;40)  \times (0;6,32,41,39,\cdots)\\ 
	&\approx  3;7,37,14,\cdots
\end{align*}
which is very close to  the approximation $\pi\approx 3;7,30$. It is also interesting that   another mathematical tablet shows that Babylonians might have also been aware of the  better approximation $\pi \approx 3;9=3.15$ (see \cite{Mur11}).

These observations suggest  that the Mesopotamian scribes  knew more accurate approximations for   $\pi $, which they did not use in everyday calculations because they could not be conveniently represented as finite sexagesimal  fractions. For a more detailed discussion  about  probable approximations of $\pi$, see \cite{Mur16}.   In what follows, beside Bruins' interpretation, we discuss two others.

\noindent
\underline{Bruins' Interpretation:}\\
He considered the constant $\frac{24}{25} =0;57,36 $ for a circle and   suggested    the following    approximate formula  for the area  $S$  of a circle  with circumference  $ c$: 
\[S\approx(0;57,36)\times \frac{c^2}{12}.\]
Since the area of a circle with circumference $c$ is
$$S=\frac{c^2}{4\pi},$$ 
it follows that 
\[\frac{c^2}{4\pi}\approx\frac{c^2}{12}\times \frac{24}{25}, \]
which implies that
\[ \pi\approx \frac{25}{8}=3.125.\]

\noindent
\underline{Neugebauer's Interpretation:}\\
 In his interpretation, he considered this number as a constant for the circumference  of a regular hexagon   and  used the approximate formula 
\[c_6\approx (0;57,36)\times c\] 
where $c_6$ and $c$ are  the circumferences   of a regular hexagon  and     its circumscribed circle  respectively. Since   
\[\frac{c_6}{c}=\dfrac{6a}{2\pi a}=\frac{3}{\pi},\]
we get
\[\frac{3}{\pi}\approx \dfrac{24}{25} \Longrightarrow \pi\approx   \frac{3\times 25}{24} =\frac{25}{8}   =3.125.\]

\begin{figure}[H]
	\centering
	\includegraphics[scale=1]{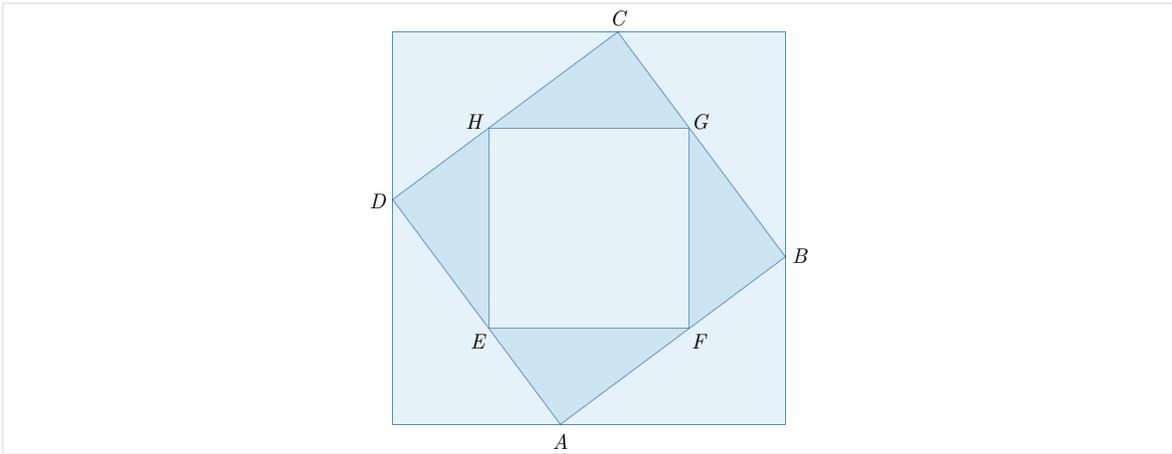}
	\caption{A chain of  right triangles}
	\label{Figure12a}
\end{figure}

\noindent 
\underline{Friberg's Interpretation:}\\
 By following \cite{Vai63}, he    offers that  this number is  the area of a chain of four right triangles  arranged in  the square  $ABCD$ with side $\overline{AB}=1 $ as shown in  \cref{Figure12a}.

\subsubsection*{Line 31: Diagonal of a Square}
\begin{figure}[H]
	\centering
	\includegraphics[scale=1]{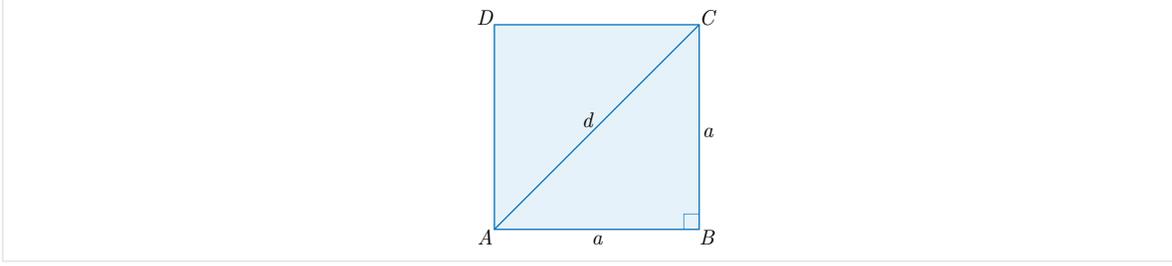}
	\caption{Diagonal of a square}
	\label{Figure12}
\end{figure}

In this line, we read {\fontfamily{qpl}\selectfont 1,25 igi-gub \textit{\v{s}\`{a}} bar-d\'{a} \textit{\v{s}\`{a}} nigin}  ``1;25 is the constant of the diagonal  of a square''.
Consider a square $ABCD$ with side $a$ and diagonal $d$ as shown in \cref{Figure12}. It follows from the Pythagorean theorem in the right triangle $\triangle ABC $  that 
\[d=\sqrt{a^2+a^2}= \sqrt{2}a .\]
These two facts implies that
\[ \sqrt{2}\approx 1;25 \]
confirming  that   the Susa scribes     knew  and used the approximate value $\sqrt{2}\approx 1;25=\frac{17}{12}$    which  is considered as one of the approximations of $\sqrt{2}$ in the Babylonian mathematics.  It should be noted that the most common approximate value for $\sqrt{2}$ in the Babylonian mathematics was $\frac{3}{2}=1;30$.

\subsubsection*{Line 32: Diagonal of a Rectangle}
\begin{figure}[H]
	\centering
	\includegraphics[scale=1]{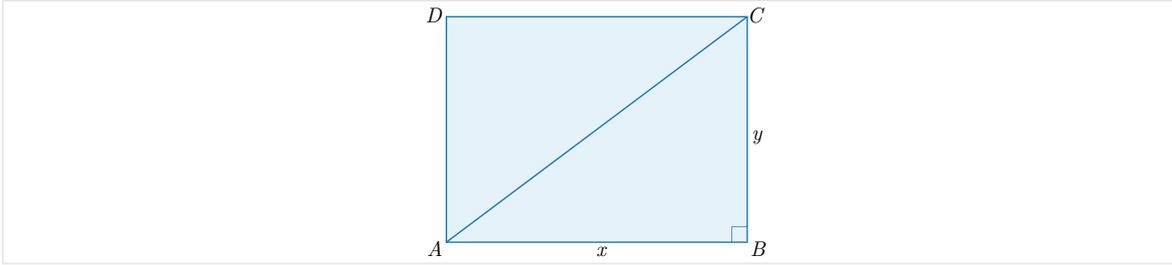}
	\caption{Diagonal of a rectangle}
	\label{Figure13}
\end{figure}

In this line, we read {\fontfamily{qpl}\selectfont 1,15 igi-gub \textit{\v{s}\`{a}} bar-d\'{a} u\v{s} \textit{\`{u}} sag} ``1;15 is the constant of the diagonal  of a rectangle''. We claim that the sides of this rectangle are of length   1 and 0;45 using  the Babylonian tradition that one side of the figure should be of length 1. In fact, if    $ \overline{AB}=1$ in \cref{Figure13}, then     by the Pythagorean theorem, we have
\[\overline{BC}=\sqrt{\overline{AC}^2-\overline{AB}^2}=\sqrt{(1;15)^2-1^2}=\sqrt{\dfrac{9}{16}}     =\dfrac{3}{4}=0;45.\]

Note that without using the assumption $ \overline{AB}=1$, we can still find   the length and the width. In fact, let $x$ and $y$ be the length and width of our rectangle respectively as well as  assume that $ \overline{AC}=(1;15)\times a=\frac{5a}{4}$, for some $a>0$. Then by the Pythagorean theorem we have $x^2+y^2=\frac{5^2a^2}{4^2} $ or equivalently
\[ \left(\frac{4x}{a}\right)^2+\left(\frac{4y}{a}\right)^2=5^2. \]
This means $(\frac{4x}{a},\frac{4y}{a},5)$  is a Pythagorean triple. So there are coprime   integers like $m>n>0$ such that they are not both odd and
\begin{equation*}
	\begin{dcases}
		\frac{4y}{a}=m^2-n^2,\\
		\frac{4x}{a}=2mn,\\
		5=m^2+n^2.
	\end{dcases}
\end{equation*}
Since $5=2^2+1^2$, it follows from the third equation that $m=2$ and $n=1$. Thus, the first two equations imply that 
$$y=\frac{(2^2-1^2)a}{4}=\frac{3a}{4}=(0;45)a $$
and 
$$x=\frac{(2\times 2\times 1)a}{4}=a$$  It should be asserted that the Pythagorean triple $(3,4,5)$ plays a key role here, because we have  
\[ 1;15=\frac{1}{4}\times 5,~~1=\frac{1}{4}\times 4,~~ \text{and}~~0;45=\frac{1}{4}\times 3.\]

\subsection{SMT No.\,15}
 So far, we have only considered the   implicit applications of the Pythagorean theorem. However,  the scribe of this tablet is explicitly using the Pythagorean rule in his calculations. The reader can   see the wording of the scribe in lines 7-10 and lines 12-15 which   exactly describe the Pythagorean rule.

\subsubsection*{Transliteration}\label{SS-TI-SMT15}

\begin{note1}
	\textbf{Obverse}\\ 
	\underline{First Problem:  Lines 1-15}\\
	(L1)\hspace{2mm} \{20\} k\'{a} 20 \textit{di-ik-\v{s}u} 2,30-ta-\`{a}m \textit{\'{u}r-t}[\textit{a-bi}]\\
	(L2)\hspace{2mm} za-e $ \frac{\text{1}}{\text{2}} $ 20 \textit{di-ik-\v{s}\'{i} he-pe} 10 \textit{ta-mar} $ \frac{\text{1}}{\text{2}} $ 30 [\textit{he-pe}]\\
	(L3)\hspace{2mm} 15 \textit{ta-mar i-na} 20 \textit{di-ik-\v{s}\'{i}} 15 zi 5 \textit{t}[\textit{a-mar}]\\
	(L4)\hspace{2mm} 10 nigin 1,40 \textit{ta-mar} igi-5 \textit{pu-\c{t}}[\textit{ur}]12 \textit{ta}-[\textit{mar} 12 \textit{a-na} 1,40]\\
	(L5)\hspace{2mm} \textit{i-\v{s}\'{i}} 20 \textit{ta-mar} $ \frac{\text{1}}{\text{2}} $ \textit{he-pe} 10 \textit{ta-mar} $ \frac{\text{1}}{\text{2}} $ [$\cdots$ $\cdots$ $\cdots$] $\cdots$ $\cdots$\\
	(L6)\hspace{2mm} 2,30 \textit{a-na} [10 dah] 12,30 \textit{ta-mar i-na} 12,30 [$\cdots$ $\cdots$ \textit{ta}]-\textit{mar}\\
	(L7)\hspace{2mm} 2,30 [\textit{i-na} 12,30] zi-\textit{ma} 10 \textit{ta-mar} 12,30 nigin 2,36,15 \textit{ta-mar}\\
	(L8)\hspace{2mm} 10 nigin 1,40 \textit{ta-mar} 1,40 \textit{i-na} 2,36,15 zi-\textit{ma}\\
	(L9)\hspace{2mm} 56,15 \textit{ta-mar mi-na} [\'{i}]b-si 7,30 \'{i}b-si\\
	(L10)\hspace{0mm} 15 \textit{ta-mar} 15 dal-1-kam dal-2-kam \textit{ki-i ta-mar}\\
	\underline{Second Problem:  Lines 11-16}\\
	(L11)\hspace{0mm} za-e 2,30 \textit{pa-na-am \v{s}\`{a}} tu-tu-da \textit{\`{u}} 2,30 2-kam ul-gar\\
	(L12)\hspace{0mm} 5 \textit{i-na} 12,30 zi-\textit{ma} 7,30 \textit{ta-mar} 12,30 nigin\\
	(L13)\hspace{0mm} 2,36,15 \textit{ta-mar} 7,30 nigin 56,15 \textit{ta-mar}\\
	(L14)\hspace{0mm} 56,15 \textit{i-na} 2,36,15 zi-\textit{ma} 1,40 \textit{ta-mar}\\
	(L15)\hspace{0mm} 1,40 \textit{mi-na} \'{i}b-si 10 íb-si 10 \textit{a-na} 2 tab-ba 20 \textit{ta-mar}\\
	(L16)\hspace{0mm} 20 dal-2-kam\\
	
	\textbf{Reverse}\\ 
	\underline{Third Problem:  Lines 1-7}\\
	(L1)\hspace{2mm} za-e $ \frac{\text{1}}{\text{2}} $ [$\cdots$ $\cdots$ $\cdots$]\\
	(L2)\hspace{2mm}  3,45 [$\cdots$ $\cdots$ $\cdots$]\\
	(L3)\hspace{2mm}  3,45 \textit{a-na} [$\cdots$ $\cdots$]\\
	(L4)\hspace{2mm}  $ \frac{\text{1}}{\text{2}} $ 11,15 \textit{he-pe} [5,37,30 \textit{ta-mar} $\cdots$]\\
	(L5)\hspace{2mm}  10 \textit{t}[\textit{a-mar}] 10 \textit{a}-[\textit{na} $\cdots$]\\
	(L6)\hspace{2mm}  15 7,30 [$\cdots$ $\cdots$]\\
	(L7)\hspace{2mm}  15 [$\cdots$ $\cdots$ $\cdots$]
\end{note1}

\subsubsection*{Translation}\label{SS-TR-SMT15}  

\noindent 
\textbf{Obverse}\\
\underline{First Problem:  Lines 1-10}
\begin{tabbing}
	\hspace{14mm} \= \kill 
	(L1)\> \tabfill{I have enlarged the gate to make 20 ({\fontfamily{qpl}\selectfont k\`{u}\v{s}}) in width by an extension of 2;30 ({\fontfamily{qpl}\selectfont k\`{u}\v{s}}) in every direction.}\index{kuz@k\`{u}\v{s} (length unit)}\index{width}\\
	(L2)\> \tabfill{You, halve 20 of the enlargement, (and) you see 10. Halve 30, (and)}\\
	(L3)\> \tabfill{you see 15. Subtract 15 from 20 of the enlargement, (and) you see 5.}\\ 
	(L4)\> \tabfill{Square 10, (and) you see 1,40. Make the reciprocal of 5, (and) you see 0;12.}\index{reciprocal of a number} \\ 
	(L5)\> \tabfill{Multiply it by 1,40, (and) you see 20. Halve 20, (and) you see 10. Halve $\cdots$}\\
	(L6)\> \tabfill{Add 2;30 to 10, (and) you see 12;30. From 12;30, $\cdots$ you see $\cdots$.}\\
	(L7)\> \tabfill{Subtract 2;30 from 12;30, and you see 10. Square 12;30, (and) you see 2,36;15.}\\
	(L8)\> \tabfill{Square 10, (and) you see 1,40. Subtract 1,40 from 2,36;15, and}\\
	(L9)\> \tabfill{you see 56;15. What is the square root? 7;30 is the square root. Multiply it by 2, (and)}\index{square root}\\
	(L10)\> \tabfill{you see 15. 15 is the first space between (that is, the original width of the gate). How do you see the second space between (that is, the width of the enlarged gate)?}\index{width}
\end{tabbing} 
\noindent 
\underline{Second Problem:  Lines 11-16}
\begin{tabbing}
	\hspace{14mm} \= \kill 
	(L11)\> \tabfill{You, add the first produced 2;30 and the second 2;30 together, (and you see 5).}\\
	(L12)\> \tabfill{Subtract 5 from 12;30, and you see 7;30. Square 12;30, (and)}\\
	(L13)\> \tabfill{you see 2,36;15. Square 7;30, (and) you see 56;15.}\\ 
	(L14)\> \tabfill{Subtract 56;15 from 2,36;15, and you see 1,40.}\\
	(L15)\> \tabfill{What is the square root of 1,40? 10 is the square root. Multiply 10 by 2, (and) you see 20.}\index{square root}\\
	(L16)\> \tabfill{20 is the second space between.}
\end{tabbing} 

\noindent 
\textbf{Reverse}\\
\underline{Third Problem:  Lines 1-7}
\begin{tabbing}
	\hspace{14mm} \= \kill 
	(L1)\> \tabfill{You, halve [$\cdots$ $\cdots$ $\cdots$]}\\
	(L2)\> \tabfill{3,45 [$\cdots$ $\cdots$ $\cdots$]}\\
	(L3)\> \tabfill{3,45 to [$\cdots$ $\cdots$ $\cdots$]}\\
	(L4)\> \tabfill{Halve 11,15, (and) you see [ 5,37;30. $\cdots$]}\\
	(L5)\> \tabfill{you see 10. 10 to [$\cdots$ $\cdots$]}\\
	(L6)\> \tabfill{15 7,30 [$\cdots$ $\cdots$]}\\
	(L7)\> \tabfill{15 [$\cdots$ $\cdots$ $\cdots$]}
\end{tabbing}

\subsubsection*{Mathematical Interpretation} 
\begin{figure}[H]
 	\centering
 	\includegraphics[scale=1]{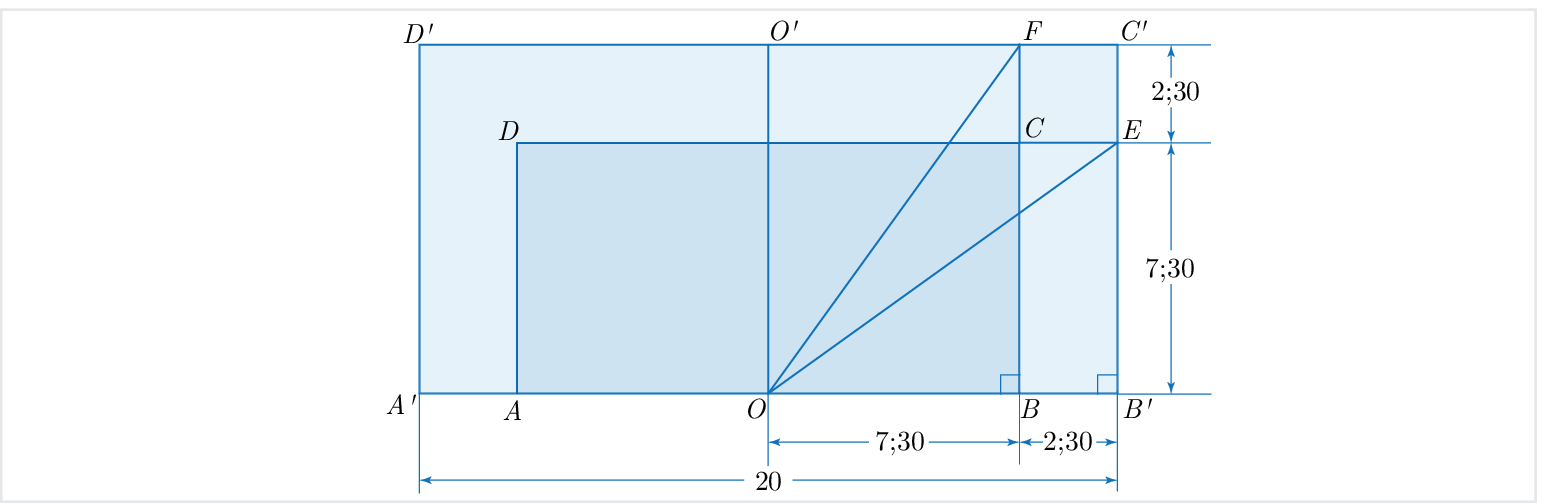}
 	\caption{Enlargement of a gate}
 	\label{Figure14}
 \end{figure}

 Not all the calculations done in this text are    clear to us. However, since we know that the subject of three problems   is the enlargement of a gate\index{enlargement of a gate}, we can  try to reconstruct the dimensions of the gate  as shown in  \cref{Figure14}.
 
The original gate (the heavily-shaded rectangle\index{rectangle}), which is 15 {\fontfamily{qpl}\selectfont k\`{u}\v{s}}\index{kuz@k\`{u}\v{s} (length unit)} ($ \approx $ 7.5m) in width\index{width} and 7;30 {\fontfamily{qpl}\selectfont k\`{u}\v{s}}\index{kuz@k\`{u}\v{s} (length unit)} in height\index{height}, has been enlarged by an extension of 2;30 {\fontfamily{qpl}\selectfont k\`{u}\v{s}}\index{kuz@k\`{u}\v{s} (length unit)} in two directions. In the first problem, the scribe of this tablet might have intended to ask for the original width\index{width} of the gate ({\fontfamily{qpl}\selectfont dal-1-kam}) when the widened width\index{width} ({\fontfamily{qpl}\selectfont dal-2-kam}) is given and the converse in the second problem, applying the Pythagorean Theorem\index{Pythagorean Theorem}  to the right triangles\index{right triangle} $\triangle OBF$ and $\triangle OB'E$ respectively (see \cref{Figure14}). But it seems that he failed to complete his calculations.

In the first problem, the first space between ({\fontfamily{qpl}\selectfont dal-1-kam}), i.e., the length\index{length} of $ AB$,  is calculated as follows. Note that the Susa  scribe\index{Susa scribes} seems to have assumed (or computed?) that $\overline{BF}=10 $, $\overline{B'E} =7;30 $ and $\overline{OF}=\overline{OE}=12;30 $. So, according to lines 7-10, we can use the Pythagorean Theorem\index{Pythagorean Theorem}  in the right triangle\index{right triangle} $ \triangle OBF$  to write
\begin{align*}
	\overline{AB}&=2\times \overline{OB}\\
	&=2\sqrt{\overline{OF}^2-\overline{BF}^2}\\
	&=2 \sqrt{(12;30)^2- 10^2}\\
	&=2 \sqrt{2,36;15-1,40}\\
	&= 2\sqrt{56;15}\\
	&= 2\sqrt{(7;30)^2}\\
	&= 2\times(7;30)\\
	&= 15.
\end{align*}

Similarly, in the second problem the second space between ({\fontfamily{qpl}\selectfont dal-2-kam}), i.e., the length\index{length} of $ A'B'$,  is obtained by using the Pythagorean Theorem\index{Pythagorean Theorem}  in the right triangle\index{right triangle} $\triangle OB'E $. According to lines 12-15,  we have
\begin{align*}
	\overline{A'B'}&=2\times \overline{OB'}\\
	&=2\sqrt{\overline{OE}^2-\overline{B'E}^2}\\
	&=2 \sqrt{(12;30)^2- (7;30)^2}\\
	&=2 \sqrt{2,36;15-56;15}\\
	&= 2\sqrt{1,40}\\
	&= 2\sqrt{(10)^2}\\
	&= 2\times 10\\
	&= 20.
\end{align*}

\subsection{SMT No.\,19}
Here, we only consider the first problem in the text which concerns an application of the Pythagorean theorem. Although the second problem uses the Pythagorean theorem in its calculations too, we   will treat it in an another article on algebraic equations. Similar to the \textbf{SMT No.\,15}, the scribe is   explicitly using the Pythagorean theorem.

\subsubsection*{Transliteration}\label{SSS-P1TI-SMT19}

\begin{note1}
	\underline{Obverse:  Lines 1-11}\\
	(L1)\hspace{2mm} sag \textit{a-na} u\v{s} \textit{re-ba-ti} $<$u\v{s}$>$ \textit{li-im-\c{t}\`{i|}} 40 \textit{\v{s}\`{a}} (text: RU) tab bar-d\'{a}\\
	(L2)\hspace{2mm} u\v{s} \textit{\`{u}} sag \textit{mi-nu} za-e 1 u\v{s} gar 1 gaba gar\\
	(L3)\hspace{2mm} 15 \textit{re-ba-ti i-na} 1 zi 45 \textit{ta-mar}\\ 
	(L4)\hspace{2mm} 1 \textit{ki-ma} u\v{s} gar 45 \textit{ki-ma} sag gar 1 u\v{s} nigin 1 \textit{ta-mar}
	\\
	(L5)\hspace{2mm} 45 sag nigin 33,45 \textit{ta-mar} 1 \textit{\`{u}} 33,45\\
	(L6)\hspace{2mm} ul-gar 1,33,45 \textit{ta-mar mi-na} \'{i}b-si 1,15 íb-si\\
	(L7)\hspace{2mm} \textit{a\v{s}-\v{s}um} 40 bar-d\'{a} \textit{qa-bu-ku} igi-1,15 bar-d\'{a} \textit{pu-\c{t}\'{u}-}$<$\textit{\'{u}r}$>$\\
	(L8)\hspace{2mm} 48 $<$\textit{ta-mar}$>$ 48 \textit{a-na} 40 bar-d\'{a} \textit{\v{s}\`{a} qa-bu-ku i-\v{s}í}
	\\
	(L9)\hspace{2mm} 32 \textit{ta-mar} 32 \textit{a-na} 1 u\v{s} \textit{\v{s}\`{a}} gar \textit{i-\v{s}\'{i}}
	\\
	(L10)\hspace{0mm} 32 \textit{ta-mar} 32 u\v{s} 32 \textit{a-na} 45 sag \textit{\v{s}\`{a}} gar\\
	(L11)\hspace{0mm} \textit{i-\v{s}í} 24 \textit{ta-mar} 24 sag
	
\end{note1}

\subsubsection*{Translation}\label{SSS-P1TR-SMT19}

\underline{Obverse:  Lines 1-11}
\begin{tabbing}
	\hspace{16mm} \= \kill 
	(L1)\> \tabfill{Let the width become less than the length by one fourth of the length. 40 is the diagonal, a partner (of the length and width).}\index{length}\index{width}\\
	(L2)\> \tabfill{What are the length and the width? You, put down 1 of the length. Put down the same number.}\index{length}\index{width}\\
	(L3)\> \tabfill{Subtract 0;15 from 1, (and) you see 0;45.}\\
	(L4)\> \tabfill{Put down 1 as the length. Put down 0;45 as the width. Square 1 of the length, (and) you see 1.}\index{length}\index{width}\\
	(L5)\> \tabfill{Square 0;45 of the width, (and) you see 0;33,45. 1 and 0;33,45,}\index{width}\\
	(L6)\> \tabfill{add (them) together, (and) you see 1;33,45. What is the square root of 1;33,45? 1;15 is the square root.}\index{square root}\\
	(L7)\> \tabfill{Since 40 of the diagonal is said to you, make the reciprocal of 1;15 of the diagonal, (and)}\index{reciprocal of a number} \\
	(L8)\> \tabfill{you see 0;48. Multiply 0;48 by 40 of the diagonal which is said to you, (and)}\\
	(L9)\> \tabfill{you see 32. Multiply 32 by 1 of the length which you put down, (and)}\index{length}\\
	(L10-11)\> \tabfill{you see 32. 32 is the length. Multiply 32 by 0;45 of the width which you put down, (and) you see 24. 24 is the width.}\index{length}\index{width} 
\end{tabbing}\index{diagonal}

\subsubsection*{Mathematical Interpretation}
 There are three variables in this problem: the length\index{length}, the width\index{width}, and the diagonal\index{diagonal}, which can be imagined as the sides and the diagonal\index{diagonal} of  a rectangle\index{rectangle}.

Let $x$, $y$ and $d$ denote the length\index{length}, the width\index{width} and the diagonal\index{diagonal} respectively. It follows from  the Pythagorean Theorem\index{Pythagorean Theorem}  that $d$  depends on $x$ and $y$:
$$d^2=x^2+y^2. $$  
Lines 1-2 give us the following system of equations:
\begin{equation}\label{equ-SMT19-aa}
	\begin{dcases}
		x-y=\dfrac{1}{4}x\\
		d= 40
	\end{dcases}
\end{equation}
or equivalently
\begin{equation}\label{equ-SMT19-a}
	\begin{dcases}
		x-y=\dfrac{1}{4}x\\
		\sqrt{x^2+y^2}=40.
	\end{dcases}
\end{equation}
 According to line 3, we can use \cref{equ-SMT19-a} to compute the value of $y$ with respect to $x$:
\begin{align*}
	&~~  x-y=\dfrac{1}{4}x \\
	\Longrightarrow~~&~~    y=x-\dfrac{1}{4}x\\
	\Longrightarrow~~&~~     y=\left(1-\dfrac{1}{4}\right)x \\
	\Longrightarrow~~&~~    y=(1-0;15)x
\end{align*}
which implies that
\begin{equation}\label{equ-SMT19-c}
	y=(0;45)x.
\end{equation} 
Next,   according to lines 5-6,  we use  \cref{equ-SMT19-a} and \cref{equ-SMT19-c} to find $d$ with respect to $x$:
\begin{align*}
	d &=\sqrt{x^2+y^2} \\
	&=\sqrt{x^2+(0;45)^2x^2}\\
	&=\sqrt{[1+(0;45)^2]x^2}\\
	&=\sqrt{(1+0;33,45)x^2}\\
	&=\sqrt{(1;33,45)x^2}\\
	&= \sqrt{(1;15)^2x^2} \\
	&=(1;15)x
\end{align*}
so
\begin{equation}\label{equ-SMT19-ca}
	d=(1;15)x.
\end{equation} 
Now,  according to lines 7-8,  the scribe finds the value of $x$. In fact, it follows from \cref{equ-SMT19-aa} and \cref{equ-SMT19-a}  that
\begin{align*}
	&~~  (1;15)x=40 \\
	\Longrightarrow~~&~~    x= \dfrac{1}{(1;15)} \times 40\\
	\Longrightarrow~~&~~     x= (0;48)\times 40
\end{align*} 
giving that
\begin{equation}\label{equ-SMT19-cc}
	x=32.
\end{equation} 
Finally, according to lines 9-11, he obtains the values of $x$ and $y$ by using  \cref{equ-SMT19-c} and \cref{equ-SMT19-cc} as follows:
\[ x= 32 \]
and
\[ y=(0;45)x=(0;45)\times 32=24. \] 
Also, the value of $d$ is easily obtained from the previous calculations:
\[ d= \sqrt{x^2+y^2}=\sqrt{17,4+9,36}=\sqrt{26,40}=40. \] 
Note that the Pythagorean triple treated in this problem is $(24,32,40)=(8\times 3, 8\times 4,8\times 5)$.

\section{Conclusion} 
 The  implicit and explicit applications of the Pythagorean theorem found  in the \textbf{SMT} show that the Elamite scribes--like their Babylonian counterparts--were fully aware of this basic theorem. They freely used the Pythagorean rule whenever their calculations involved  computing a side of a right triangle.  
 
Besides the algebraic applications of this theorem in the \textbf{SMT}, the scribe of \textbf{SMT No.\,1} has presented a geometric application of the theorem. Although the scribe has only given the numerical data on the tablet, its mathematical interpretation  clearly confirms that obtaining these numbers requires the application of the Pythagorean theorem.  In fact, this text and  the  Babylonian text \textbf{YBC 7289} might be the only  geometric applications of the Pythagorean theorem in Elamite and Babylonian mathematics.

 \section*{Appendix: Mathematical Tablet BM 85196, No.\,9}\phantomsection\label{appendix}
 In this appendix, we give the transliteration, the translation and the mathematical interpretation  the 9th problem in lines 7-16  of  \textbf{BM 85196}.
 
 \subsection*{Transliteration}

 \begin{note1} 
 	\underline{Lines 7-16}\\
 	(L7)\hspace{4mm}    giš \textit{pa-lu-um} 30 gi \textit{i-na i}-[\textit{ga-ri-im} ur]-bi \textit{š}[\textit{a-ki-in}]\\
 	(L8)\hspace{4mm}    \textit{e-le-nu} 6 \textit{ur-dam} \textit{i-na ša-a}[\textit{p-la-n}]\textit{u}-[\textit{um} en-nam \textit{is-sé-a-am}] \\
 	(L9)\hspace{4mm}     za-e 30 nigin 15 \textit{ta-mar} 6 \textit{i-n}[\textit{a}] 30 ba-[zi 24 \textit{ta-mar}]\\
 	(L10)\hspace{2mm}   24 nigin 9,36 \textit{ta-mar} 9,[36 \textit{i-na} 15 ba-zi]\\
 	(L11)\hspace{2mm}    5,24 \textit{ta-mar} 5,24 en-nam [íb-si$ _{8} $ 18 íb-si$ _{8} $ 18]\\
 	(L12)\hspace{2mm}  \textit{i-na qá-qá-ri is-sé-a-am šum-ma} 18 \textit{i-n}[\textit{a qá}]-\textit{qá-ri-im} \\
 	(L13)\hspace{2mm}  \textit{e-le-nu-um} en-nam \textit{ur-dam} 18 nigin 5,24 \textit{ta-mar} \\
 	(L14)\hspace{2mm}  5,24 \textit{i-na} 15 ba-zi 9,36 \textit{ta-mar} 9,36\\
 	(L15)\hspace{2mm}   en-nam íb-si$ _{8} $ 24 íb-si$ _{8} $ 24 \textit{i-na} 30 ba-zi\\ 
 	(L16)\hspace{2mm} 6 \textit{ta-mar ur-dam ki-a-am ne-pé-šu}m  
 \end{note1}

 \subsection*{Translation}

 \underline{Lines 7-16}
 \begin{tabbing}
 	\hspace{18mm} \= \kill 
 	(L7)\> \tabfill{A timber. (Its length is) 0;30 ({\fontfamily{qpl}\selectfont nindan}, that is, 1) {\fontfamily{qpl}\selectfont gi}. At a wall it is placed vertically.}\\
 	(L8)\> \tabfill{(From) above I went down by 0;6 ({\fontfamily{qpl}\selectfont nindan}). How far is it (the lower end of the timber) from the base (of the wall)?} \\ 
 	(L9)\> \tabfill{You, square 0;30, and you see 0;15. Subtract 0;6 from 0;30, and you see 0;24.}\\ 
 	(L10)\> \tabfill{Square 0;24, and you see 0;9,36. Subtract 0;9,36 from 0;15, and}\\
 	(L11)\> \tabfill{you see 0;5,24. What is the square root of 0;5,24? 0;18 is the square root.}\\ 
 	(L12)\> \tabfill{It is 0;18 ({\fontfamily{qpl}\selectfont nindan}) away from the bottom. If it is 0;18 ({\fontfamily{qpl}\selectfont nindan}) away from the bottom,}\\
 	(L13)\> \tabfill{how far did I go down from above? Square 0;18, and you see 0;5,24.}\\
 	(L14)\> \tabfill{Subtract 0;5,24 from 0;15, and you see 0;9,36.}\\ 
 	(L15)\> \tabfill{What is the square root of 0;9,36? 0;24 is the square root. Subtract 0;24 from 0;30,}\\
 	(L16)\> \tabfill{and you see 0;6. I went down (by 0;6 from above). Such is the procedure.}
 \end{tabbing}

 \subsection*{Mathematical Interpretation}
 In this text, the scribe is dealing with a situation pictured in \cref{Figure4}. First, he assumes $h_0=0;6$, $l=h=0;30$ and computes $d$ as
 \begin{align*}
 	d&=\sqrt{l^2-(h-h_0)^2}\\
 	 &=\sqrt{(0;30)^2-(0;30-0;6)^2}\\
 	 &=\sqrt{(0;30)^2-(0;24)^2}\\
 	 &=\sqrt{0;15-0;9,36}\\
 	 &=\sqrt{0;5,24}\\
 	 &=0;18.
 \end{align*}

Then, he assumes $l=h=0;30$, $d=0;18$ and computes $h_0$: 
  \begin{align*}
 	h_0 &=h-\sqrt{l^2-d^2}\\
 	&=0;30-\sqrt{(0;30)^2-(0;18)^2}\\
 	&=0;30-\sqrt{0;15-0;5,24}\\
 	&=0;30-\sqrt{0;9,36}\\
 	&=0;30-0;24\\
 	&=0;6.
 \end{align*}

\end{document}